\magnification1200
\def\label#1/{\ifmmode\leqno(#1)\else#1\fi}
\def\cite#1{\excite#1::/}
\def\excite#1:#2:#3/{\edef\griff{cit#1}%
\formcite{\plazieres\def\next{#2}\ifx\next\empty\else\formciteaddl{#2}\fi}}
\def\formcite#1{{\bf[}#1{\bf]}}
\def\formciteaddl#1{: #1}
\def\plazieres{\expandafter\ifx\csname\griff\endcsname\relax%
\xdef\esfehlt{\griff}\blackmark\else{\csname\griff\endcsname}\fi}
\def\nextref{\item{[\rh]}}
\def\blackmark{\immediate\write16{undefined forward reference: \esfehlt}}

\def \lemmicchelli{1}
\def \remmicchelli{2}
\def \basicidentity{3}
\def \spannings{4}
\def \propone{5}
\def \corone{6}
\def \corip{7}
\def \lemschaback{8}
\def \grammatrix{9}
\def \changeofvar{10}
\def \neededcond{11}
\expandafter\def \csname citBR1\endcsname{BR1}
\expandafter\def \csname citBR2\endcsname{BR2}
\expandafter\def \csname citM\endcsname{M}
\expandafter\def \csname citSX\endcsname{SX}
\expandafter\def \csname citS\endcsname{S}

\def\NN{{{\rm I}\kern-.16em {\rm N}}}
\def\RR{{{\rm I}\kern-.16em {\rm R}}}
\def\ZZ{{{\rm Z}\kern-.28em{\rm Z}}}
\def\FF{{{\rm I}\kern-.16em {\rm F}}}
\def\Rd{\RR^d}
\def\all{\mathord{\rm all}\;}
\let\dword\defw
\let\eword\empw
\def\formal#1{\bigskip{\bf #1}\hskip1em}
\def\Proof#1{\indent{\bf Proof#1:}\qquad}
\def\proof{\Proof{}}
\def\makeblanksquare#1#2{
\dimen0=#1pt\advance\dimen0 by -#2pt
      \vrule height#1pt width#2pt depth0pt\kern-#2pt
      \vrule height#1pt width#1pt depth-\dimen0 \kern-#1pt
      \vrule height#2pt width#1pt depth0pt \kern-#2pt 
      \vrule height#1pt width#2pt depth0pt
}

\let\endproofsymbol\boxx
\def\endproof{\hbox{}~\hfill\endproofsymbol\medskip\goodbreak}
\let\gL\Lambda
\def\gM{{\rm M}}
\let\ga\alpha
\let\gb\beta
\let\gk\kappa
\let\gl\lambda
\def\heading#1{\bigskip\goodbreak\centerline{\bf #1}
               \nobreak\medskip\nobreak}
\def\inpro#1{\langle#1\rangle}
\def\inv#1{#1^{-1}}
\def\newpar{\par}
\def\norm#1{\Vert#1\Vert}
\def\unitv#1{{\bf i}_{#1}}
\def\ran{\mathop{\rm ran}\nolimits}\def\Zp{\ZZ_+}
\def\Zdp{\Zp^d}

\let\Span\spam
\def\braket#1{\hbox{$[\![$}#1\hbox{$]\!]$}}
\def\tp{{}^{\rm t}}
\def\title{On interpolation by radial polynomials}
\def\author{C. de Boor}
\centerline{\bf\title}
\smallskip
\centerline{\author}
\smallskip
\bigskip
\formal{Abstract}
A lemma of Micchelli's, concerning radial polynomials and weighted sums of 
point evaluations, is shown to hold for arbitrary 
linear functionals, as is Schaback's more recent extension of this
lemma and Schaback's result concerning interpolation by radial polynomials.
Schaback's interpolant is explored.

\bigskip
In his most-cited paper, \cite{M}, Micchelli supplies the following 
interesting auxiliary lemma (his Lemma 3.1).

\proclaim (\label\lemmicchelli/) Lemma. If $\sum_{i=1}^n c_i p(x_i) = 0$ for all
$p\in\Pi_{<k}(\Rd)$, then
$$(-1)^k\sum_{i=1}^n\sum_{j=1}^n c_ic_j\norm{x_i-x_j}^{2k} \ge 0,$$
where equality holds if and only if
$$
\sum_{i=1}^n c_ip(x_i)=0,\quad p\in \Pi_{\le k}(\Rd).$$

Here, $(x_1,\ldots,x_n)$ is a sequence in $\Rd$, and $\norm{\cdot}$ is the
Euclidean norm in $\Rd$. Further, with
$$()^\ga:\Rd\to\RR:x\mapsto x^\ga:=x(1)^{\ga(1)}\cdots x(d)^{\ga(d)}$$
a convenient if nonstandard notation for the power function,
$$\Pi_{\le k} := \Pi_{\le k}(\Rd):=\Span(()^\ga: \ga\in\Zdp, |\ga|\le k)$$
is the collection of
all polynomials in $d$ real variables of (total) degree $\le k$, and
$$\Pi_{<k} := \Pi_{<k}(\Rd)$$
those of degree $<k$. The collection of all polynomials in
$d$ real variables will be denoted here, correspondingly, by $\Pi = \Pi(\Rd)$.

Micchelli follows the lemma by the following

\proclaim (\label\remmicchelli/) Remark. Applying the Lemma inductively 
shows that the conditions
$$
\sum_{i=1}^n c_ip(x_i)=0,\quad p\in \Pi_{\le k}(\Rd),$$
and
$$(-1)^k\sum_{i=1}^n\sum_{j=1}^n c_ic_jq(\norm{x_i-x_j}^{2k})=0,\quad
q\in\Pi_{\le k}(\RR),$$
are equivalent.

The essence of the proof is, perhaps, the observation
(implicit in Micchelli's proof) that
$$\eqalign{\norm{x-y}^{2k}
&\;=\;
(\norm{x}^2 - 2\sum_ix(i)y(i) + \norm{y}^2)^k\cr
&\;=\;
\sum_{a+|\gb|+c=k}{k!\over a!\gb!c!}\norm{x}^{2a}(-2)^{|\gb|}x^\gb y^\gb\norm{y}^{2c}\cr
&\;=\;
\sum_{a+b+c=k}(-2)^b\sum_{|\gb|=b}{k!\over a!\gb!c!}p_{a,\gb}(x)p_{c,\gb}(y),\cr
}\label\basicidentity/$$
with
$$p_{a,\gb}(x):=\norm{x}^{2a}x^\gb$$
and
with $a,b,c$ nonnegative integers, $\gb\in\Zdp$, and 
$$|\gb|:=\sum_j\gb(j),\quad \gb!:=\gb(1)!\cdots\gb(d)!.$$
Each summand in the final sum of (\basicidentity) is thus the product of a 
(homogeneous) polynomial in $x$ of degree $2a+b$ and a (homogeneous) polynomial 
in $y$ of degree $2c+b$.  Hence, if 
$$\gl\perp\Pi_{<k},$$
i.e.,  $\gl$ is a linear functional on $\Pi$ that vanishes on
$\Pi_{<k}$, then the tensor product of $\gl$ with itself, i.e., the linear map 
$$\gl\otimes\gl:\Pi\otimes\Pi\to\RR:()^\ga\otimes()^\gb\mapsto
\gl()^\ga\,\gl()^\gb,$$
annihilates all the
summands with $2a+b<k$ \eword{or} $2c+b<k$. As to any other summands, they must have
$2a+b\ge k$ \eword{and} $2c+b\ge k$, hence 
$$2k = 2(a+b+c) = 2a+b + 2c+b \ge 2k,$$
therefore 
$$2a+b=k=2c+b.$$
These are the summands in which the polynomial in $x$ equals the polynomial
in $y$ and, moreover, $k-b$ is even. Thus, altogether (and in faulty but
understandable notation),
$$(\gl\otimes\gl)\norm{x-y}^{2k} = 
(-1)^k\sum_{k-b\ \rm even}2^b{k!\over (((k-b)/2)!)^2}\sum_{|\gb|=b}(\gl p_{(k-b)/2,\gb})^2/\gb!.
$$
In particular, $(\gl\otimes\gl)((-1)^k\norm{x-y}^{2k})\ge0$ with equality if
and only if $\gl$ vanishes on the sequence
$$(p_{(k-b)/2,\gb}: |\gb|=b; b\in\Zp,\; k-b\ \ \rm even).
\label\spannings/
$$
But since each $p_{(k-b)/2,\gb}$ here is homogeneous of degree $k$ while,
in particular, $k-b$ is even when $b=k$ hence each $()^\gb$ with $|\gb|=k$
appears in (\spannings), this last condition is
equivalent to having $\gl$ vanish on $\Span(()^\gb: |\gb|=k,
\gb\in\Zdp)$, hence on $\Pi_{\le k}$.

Altogether, this proves the following generalization of Micchelli's Lemma and
Remark.

\proclaim (\label\propone/) Proposition.
If $\gl\perp\Pi_{<k}$, then
$$
(-1)^k(\gl\otimes\gl)\norm{x-y}^{2k}\ge 0
$$
with equality iff $\gl\perp\Pi_{\le k}$.

\proclaim (\label\corone/) Corollary.
$\gl\perp\Pi_{\le k}$ if and only if
$$(\gl\otimes\gl)\norm{x-y}^{2r}=0,\quad 0\le r\le k.$$

Note the following useful

\proclaim (\label\corip/) Corollary. The bilinear form
$$\inpro{,}_k:\Pi'\otimes\Pi':
(\gl,\mu)\mapsto(-1)^k(\gl\otimes\mu)\norm{x-y}^{2k}$$
on the algebraic dual $\Pi'$ of $\Pi$
is an inner product on any algebraic complement of 
$${\perp}\Pi_{\le k}:=\{\mu\in\Pi': \mu\perp\Pi_{\le k}\}$$
in ${\perp}\Pi_{<k}$.

Schaback \cite{S} reiterates Micchelli's results and extends them as follows
(Lemmata 8 and 9 of \cite{S}, though only for linear functionals $\gl$ that are
linear combinations of point evaluations).

\proclaim (\label\lemschaback/) Lemma.
If $\gl\perp \Pi_{\le k}$, then, for all $k\le 2\ell$, 
$$x\mapsto \gl\norm{x-\cdot}^{2\ell}$$
has degree $< 2\ell-k$.
\newpar
Conversely, if, for some $k\le2\ell$,
$$x\mapsto \gl\norm{x-\cdot}^{2\ell}$$
has degree $<2\ell-k$, then 
$\gl\perp \Pi_{\le k}$.

To be sure, the first assertion follows directly from the basic identity
(\basicidentity) since, by that identity, application of such $\gl$ to
$\norm{x-\cdot}^{2\ell}$ kills all summands with $2c+|\gb|\le k$, leaving only
those with $2\ell - 2a - |\gb| = 2c+|\gb|> k$, i.e., with $2\ell-k>
2a+|\gb|$.

For the second assertion, simply ``apply the
same idea as in the proof of Micchelli's lemma'', to quote \cite{S}. 
Arguing perhaps differently, rewrite (\basicidentity) in terms of polynomial 
degrees in $x$  to get
$$\norm{x-y}^{2\ell} = \sum_{j=0}^{2\ell}\sum_{2a+b=j}(-2)^b
\sum_{|\gb|=b}{\ell!\over a!\gb!c!}p_{a,\gb}(x)p_{c,\gb}(y).
$$
Now, $(-2)^b = (-1)^j2^b$ since $j-b$ here is always even.
Also (using $a+b+c=\ell$),
$p_{a,\gb} = \norm{\cdot}^{2(j-\ell)}p_{c,\gb}$. Therefore,
$$\norm{x-y}^{2\ell} =
\sum_{j=0}^{2\ell}(-1)^j\norm{x}^{2(j-\ell)}\sum_{2a+b=j}2^b
\sum_{|\gb|=b}{\ell!\over a!\gb!c!}p_{c,\gb}(x)p_{c,\gb}(y).
$$
Hence, 
if now $x\mapsto\gl\norm{x-\cdot}^{2\ell}$ is of degree $<2\ell-k$, then
each of the sums
$$ \sum_{2a+b=j}2^b
\sum_{|\gb|=b}{\ell!\over a!\gb!c!}p_{c,\gb}\gl
p_{c,\gb}, \quad 2\ell-k\le j,$$
must be zero, hence so must be the value of $\gl$ on each such
sum, i.e., 
$$0\;=\; \sum_{2a+b=j}2^b
\sum_{|\gb|=b}{\ell!\over a!\gb!c!}(\gl p_{c,\gb})^2,
\quad 2\ell-k\le j.$$
This implies that
$$\gl p_{c,\gb}=0,\quad |\gb|=b,\; 2\ell-k\le 2a+b,\; \ell= a+b+c,\;
0\le a,b,c,$$
hence, for the particular choice $\ell=a+b$,
therefore $c=0$ and $2\ell-2a-b\le k$, i.e., $b\le k$,
$$\gl p_{0,\gb} = 0, \quad |\gb|\le k.
$$
But that says that $\gl\perp\Pi_{\le k}$.

\bigskip In what follows, it is convenient to consider, for any sequence
or indexed `set' $(v_j: j\in J)$ of vectors in some linear space $Y$ over the
scalar field $\FF$, the corresponding map
$$[v_j: j\in J]:\FF_0^J\to Y: c\mapsto \sum_{j\in J} c(j)v_j,$$
with $\FF_0^J$ denoting all scalar-valued functions on $J$ with finite support. 
Note that the $v_j$ enter the description of this linear map in
exactly the manner in which the columns of an $(m,n$)-matrix $A$ enter the
description of the corresponding linear map $\FF^n\to\FF^m:c\mapsto Ac$,
hence it seems reasonable to call $[v_j: j\in J]$ the \dword{column map with
columns $v_j$}. We can think of $[v_j: j\in J]$ as a row, much as we can
think of a matrix as the row of its columns.
The sequence $(v_j: j\in J)$ is a basis for $Y$ exactly when $[v_j: j\in J]$
is invertible, in which case one might just as well refer to the latter as a 
basis for $Y$.
Further, if $(\gl_i: i\in I)$ is an indexed `set' in the dual, $Y'$, of $Y$,
then it is convenient to consider the corresponding map
$$[\gl_i: i\in I]\tp :Y\to\FF^I: y\mapsto(\gl_i y: i\in I),$$
calling it the \dword{row map with rows $\gl_i$} (hence the use of the
transpose sign) since the $\gl_i$ enter the description  of this linear map
in exactly the same manner in which the rows of the transpose $A\tp$ (i.e., the columns of $A$) of an $(m,n)$-matrix $A$ 
 enter the description of the corresponding linear
map $\FF^m\to\FF^n:c\mapsto A\tp c$. We can think of $[\gl_i: i\in I]\tp$ as
a column, much as we can think of a matrix as the column of its rows.
With this, we are ready to think, as we may, of
the composition of such a row map with such a column map as matrix, i.e., the
\dword{Gram matrix} or \dword{Gramian} 
$$[\gl_i: i\in I]\tp[v_j: j\in J] = (\gl_iv_j: i\in I, j\in J).$$
This use of the superscript ${}\tp$ also seems consistent 
with the notation $x\tp y:=\sum_{j=1}^d x(j)y(j)$ for the 
scalar product of $x,y\in\Rd$, given that it is standard to think of the
elements of $\Rd$ as columns.
For completeness, we note that the composition $[v_j: j\in J][\gl_i: i\in
I]\tp$ makes sense only when $I=J$, in which case it is a linear map from $Y$
to $Y$ and the most general such in case $Y$ is finite-dimensional.
\bigskip
Schaback \cite{S} considers, for the $n$-dimensional space $\gM$ spanned by
evaluation at the elements of the given $n$-set $X$ in $\Rd$, a basis
$\gL=[\gl_1,\ldots,\gl_n]$
\dword{graded (by degree)} (he calls any such a `discrete moment basis') in the sense
that the sequence 
$(\gk_i:\,i=1,\ldots,n\,)$,
with
$$\gk_i:=\max\{k: \gl_i\perp \Pi_{<k}\},\quad\all i,$$
is nondecreasing, and, for each $k$, $[\gl_i: \gk_i\ge k]$ is a basis for 
$\gM\cap{\perp}\Pi_{<k}$.

One readily obtains such a basis from any particular basis
$[\mu_1,\ldots,\mu_n]$ for $\gM$,
by applying Gauss elimination with row interchanges to the
Gram matrix
$$(\mu_i()^\ga: i=1,\ldots,n, \ga\in\Zdp) = [\mu_1,\ldots,\mu_n]\tp  V,\label\grammatrix/$$
with the columns of
$$V:=[()^\ga: \ga\in\Zdp]
\;=:\;[()^{\ga_j}: j=1,2,\ldots]$$
so ordered that $j\mapsto |\ga_j|$
is nondecreasing.

Indeed, Gauss elimination applied to an onto Gram matrix such as 
(\grammatrix) can be interpreted as providing an invertible matrix $L$ (the
product of a permutation matrix with a lower triangular matrix) and
thereby the basis
$$[\gl_1,\ldots,\gl_n] := [\mu_1,\ldots,\mu_n](\inv{L})\tp $$
for $\gM$,
and a subsequence $(\gb_1,\ldots,\gb_n)$ of $(\ga_j: j=1,2,\ldots)$
so that,
for each $i$, the first nonzero entry in row $i$ of 
$[\gl_1,\ldots,\gl_n]\tp V$, i.e., in $\gl_iV$, is the $\gb_i$th. 
Since, by assumption, the map
$j\mapsto|\ga_j|$ is nondecreasing, this implies that
$$\gk_i = \max\{k: \gl_i\perp\Pi_{<k}\} = |\gb_i| ,\quad i=1,\ldots,n,$$
and, in particular, $i\mapsto\gk_i$ is nondecreasing. Further,
if $\sum_ic(i)\gl_i \perp\Pi_{<k}$, then, since $\gk_i=|\gb_i|$, also 
$\sum_{|\gb_i|<k}c(i)\gl_i\perp\Pi_{<k}$. But this implies that
$$[c(i): |\gb_i|<k]B = 0,$$
with the matrix
$$B:=(\gl_i()^{\gb_j}: |\gb_i|,|\gb_j|<k)$$
square upper triangular with nonzero diagonal entries, hence invertible, and
therefore $c(i)=0$ for $|\gb_i|<k$. This proves that,  for each $k$, 
$[\gl_i: |\gb_i|\ge k]$ is a basis for $\gM\cap{\perp}\Pi_{<k}$.

Schaback then considers the polynomials
$$w_j:x \mapsto \gl_j\norm{x-\cdot}^{2\gk_j},\quad j=1,\ldots,n.
$$
Note that, by (\lemschaback)Lemma, 
$$\deg w_j = 2\gk_j - \gk_j = \gk_j.$$
This implies that the Gram matrix
$$\gL\tp  W=(\gl_i w_j: i,j=1,\ldots,n)$$
is block upper triangular since
$$(\gL\tp  W)(i,j)=\gl_i w_j =
(\gl_i\otimes\gl_j)\norm{x-y}^{2\gk_j}=(-1)^{\gk_j}\inpro{\gl_i,\gl_j}_{\gk_j}$$
is zero as soon as $\gk_i>\gk_j$. Further, for each $k$ and with
$$I_k:=\{i: \gk_i=k\},$$
the diagonal block
$$\gL\tp W(I_k,I_k) = (-1)^k(\inpro{\gl_i,\gl_j}_k: i,j\in I_k)$$
is invertible, by (\corip)Corollary and the fact that $[\gl_i: i\in I_k]$ is
a basis for an algebraic complement of ${\perp}\Pi_{\le k}$ in 
${\perp}\Pi_{<k}$.  We will use later that
this implies that $W$ is a \dword{graded} basis for 
$$F:=\ran W := \{\sum_j a(j)w_j: a\in\RR^n\}$$
in the sense that $j\mapsto \deg w_j$ is nondecreasing and, for each $k$,
$[w_j: \deg w_j<k]$ is a basis for $F\cap\Pi_{<k}$.

For the moment, we only use the conclusion (which Schaback draws
for the case that $\gM$ is spanned by point evaluations) that $\gL\tp  W$ is 
invertible, hence
$$P = P_S:= W\inv{(\gL\tp  W)}\gL\tp $$
is the linear projector that associates with each $p\in\Pi$ the unique element
$f\in F$ that \dword{matches $p$ at $\gM$} in the sense that 
$$\mu f=\mu p,\quad \mu\in\gM.$$

Schaback also observes that $P_S$ is of minimal degree in the sense
that $F$ minimizes
$$\deg G := \max\{\deg g: g\in G\}$$
among all polynomial subspaces $G$ that are \dword{correct for} $\gM$ in the
sense that, for every $p\in\Pi$, they contain
a unique match at $\gM$.

However, much more is true. $P_S$
is of minimal degree in the strong sense (of, e.g., \cite{BR2}) that it is  \dword{degree-reducing}, meaning that
$$\deg P_S p\le\deg p,\quad p\in\Pi.$$
This condition is shown, in \cite{BR2}, to be equivalent to the
following property more explicitly associated with the words `minimal degree':

\proclaim Definition. The finite-rank linear projector $P$ on $\Pi$
is \dword{of minimal degree}$\;\;:=\;\;$
$$\dim(G\cap\Pi_{<k})\le \dim(\ran P\cap\Pi_{<k}),\quad k\in\NN,$$
for all linear subspaces $G$ of $\Pi$ that are
correct for ${\perp}\ker P$.

\proclaim Proposition. $P_S$ is of minimal degree.

\proof Let $G$ be a linear subspace of $\Pi$ correct for 
$\gM$.
Then, $G$ is $n$-dimensional and, for any bases $\gL$ of
$\gM$ and $W$ of $G$, respectively, the Gramian $\gL\tp W$ is 
invertible.

Choose, in particular, $\gL$ to be a graded basis for 
$\gM$ and $W$ to be a graded basis for $G$. Then, for any 
$k\in\NN$, the first $\dim(G\cap\Pi_{<k}) = \#\{j: \deg w_j<k\}$ columns of the 
Gramian $\gL\tp W=(\gl_iw_j)$ have nonzero entries only in the first $\#\{i: 
\gk_i<k\} = n-\#\{i:\gk_i\ge k\} = n - \dim(\gM\cap{\perp}\Pi_{<k})$ rows. The
invertibility of the Gramian therefore implies that
$$\dim(G\cap\Pi_{<k}) \le  n - \dim(\gM\cap{\perp}\Pi_{<k}).$$
On the other hand, there is equality here when $G=\ran P_S$ since, as observed
earlier, Schaback's $w_j$ form a
graded basis while $\gk_i=\deg w_i$, all $i$, hence
$n-\dim(\gM\cap{\perp}\Pi_{<k}) = \#\{i: \gk_i<k\} =\#\{j:\deg w_j<k\}$.
\endproof

The proof shows that a linear projector $P$ on $\Pi$ is of minimal degree 
if and only if
$$\dim (\ran P\cap\Pi_{<k})\;+\;\dim({\perp}\ker P\cap{\perp}\Pi_{<k}) =
\dim\ran P,\quad \all k.$$

The polynomial interpolant $P_Sp$ to $p$ at $\gM$ is,
in general (see below), not the least interpolant $P_{BR}p$ of \cite{BR2} to $p$
at $\gM$, hence we are free to give it a name, and \dword{Schaback interpolant} 
seems entirely appropriate (hence the suffix $S$).  

We now compare $P_S$ and $P_{BR}$ in the specific context of $\cite{S}$, i.e.,
when $\gM$ is spanned by evaluation on some $n$-set $X$ in $\Rd$.
In addition to being of minimal degree,
each interpolant is constant in any direction perpendicular to the affine
hull
$$\flat X$$
of $X$ (or, \dword{flat} spanned by $X$), i.e., both satisfy
$$P f(x) = (Pf)(P_X x),$$
with $P_X$ the orthoprojector of $\Rd$ onto $\flat X$.

For the Schaback interpolant, this follows from the fact that, for any $x$, and
any $y\in\flat X$,
$$\norm{x-y}^2 = \norm{P_Xx-y}^2 + \norm{x-P_Xx}^2,$$
hence, since $\gl_j\perp\Pi_{<\gk_j}$,
$$w_j(x) = \gl_j(\norm{P_Xx-\cdot}^2 + \norm{x-P_Xx}^2)^{\gk_j} 
= \gl_j\norm{P_Xx-\cdot}^{2\gk_j},$$
by (\lemschaback)Lemma.

This readily implies that {\sl both interpolants coincide in case $X$ is 
contained in a 1-dimensional flat}.

Also, {\sl both projectors commute with 
translation}, i.e., 
$$P p(\cdot+y) = (P p)(\cdot+y),$$
{\sl and interact with any unitary change of variables as follows:
$$P p(A\cdot) = (P p)(A\tp \cdot)\label\changeofvar/$$
for all real unitary matrices 
$A$}, as follows for Schaback's projector directly from the 
observation that, for all such $A$ and any $x$ and $y$,
$$\norm{x-Ay}^2 = \norm{A(A\tp x-y)}^2 = \norm{A\tp x-y}^2.$$
However, while (\changeofvar) holds for $P=P_{BR}$ and arbitrary invertible
$A$ (see, e.g., \cite{BR1}), this is not in general so for $P=P_S$,
due to the fact that it is the `kernel' $(x,y)\mapsto\exp(x\tp y)$ on which
$P_{BR}$ is based rather than Schaback's $(x,y)\mapsto\norm{x-y}^{2\ell}$).
We would therefore expect the two interpolants in general to differ when 
$\dim\flat X>1$ and $\ran P$ isn't just some $\Pi_{<k}$.

The simplest such example occurs when $X=\{x_1,\ldots,x_4\}$ 
is a $4$-set in $\RR^2$ that spans 
$\RR^2$. In that case, the range of each projector is of the form
$$\Pi_{<2}+w\RR$$
for some homogeneous quadratic polynomial $w$. For the least interpolant, 
\cite{BR1} gives $w$ as
$$x\mapsto \gl(x\tp \,\cdot\,)^2,$$
with
$$\gl: f\mapsto \sum_{j=1}^4 a(j)f(x_j)$$
such that
$$\sum a(j)x_j = 0,\quad \sum_j a(j) = 0.\label\neededcond/$$
Since $X$ spans $\RR^2$, this says that $\gl\perp\Pi_{<2}$
and, since $\#X=4$, $\gl$ is, up to scalar
multiples, the unique such element of $\gM$.
That means that Schaback's $w$ is the leading term of
$$x\mapsto\gl\norm{x-\cdot}^4.$$
Since  $\gl\perp\Pi_{<2}$, and
 $\norm{x-x_j}^4 = (\norm{x}^2 -2 x\tp x_j + \norm{x_j}^2)^2$,
this leading term is
$$x\mapsto\gl(2\norm{x}^2\norm{\cdot}^2 + 4(x\tp \,\cdot\,)^2),$$
and this is, offhand, not just a scalar multiple of the least's $w$.
For, there is no reason to believe that $\gl\norm{\cdot}^2 = 0$.

E.g., with the specific choice
$$X = (0, \unitv1, \unitv2, z)$$
involving $\unitv1:=(1,0)$ and $\unitv2:=(0,1)$,
we have 
$$\gl: f \mapsto f(z) - z(1)f(\unitv1) - z(2)f(\unitv2),$$
hence
$$\gl\norm{\cdot}^2 = \norm{z}^2 - z(1) - z(2) = z(1)(z(1)-1)
+z(2)(z(2)-1),$$
and this is zero only for special choices of $z$. To be sure, it is zero when
$z=\unitv1+\unitv2$, i.e., in case of gridded data, giving us then bilinear
interpolation.

To be sure, graded bases of the space $\gM$, of linear functionals at which to
match given values by polynomials, have been used before in multivariate 
polynomial interpolation. For example, the multivariate `finite differences' 
introduced and used in \cite{SX} are easily seen to form such a graded basis.

In particular, the construction of the least interpolant makes use of a
graded basis $\gL$ for $\gM$ (constructed by a more stable variant of Gauss 
elimination, namely Gauss elimination `by segments') but obtains the polynomial 
space $G$ of interpolants as the span of the `leasts' of the $\gl_i$. To 
recall (e.g., from \cite{BR2}), any $\gl\in\Pi'$ is uniquely representable, 
with respect to the bilinear form
$$\RR\braket{x}\otimes\Pi\to\RR:(f,p)\mapsto
\sum_{\ga\in\Zdp}D^\ga(f)(0)D^\ga p(0)/\ga!,$$
by the formal power series
$$\hat\gl:=\sum_{k=0}^\infty
\hat\gl^{[k]},$$
with
$$\hat\gl^{[k]}:=\sum_{|\ga|=k}(\gl()^\ga)()^\ga/\ga!,\quad \all k.$$
Then
$$\gk := \max\{k: \gl\perp\Pi_{<k}\} = \min\{k: \hat\gl^{[k]}\not=0\}$$
is known as the \dword{order} of $\gl$,
and $\hat\gl^{[\gk]}$ is, by definition, the \dword{least} of $\gl$ (with $\gk$
taken to be $-1$ when $\gl=0$ and, correspondingly, $\hat\gl^{[-1]}=0$).
It is easy to see that $G=\Span\{\hat\gl_i^{[\gk_i]}: i=1,\ldots,n\}$ depends 
only on $\gM$ and not on the particular graded basis $\gL$ for $\gM$ used, and 
is spanned by homogeneous polynomials, hence, equivalently, is
dilation-invariant. For the most striking properties of $G$ 
(such as a finite list of constant coefficient differential operators whose 
joint kernel is $G$), see \cite{BR1} or \cite{BR2}.

\formal{Acknowledgements} Thanks are due to Tomas Sauer for a constructive
reading of what I thought was the final draft. Further, on receiving a preprint
of the present note, Robert Schaback informed me that, in the meantime, he,
too, had extended (\propone)--(\lemschaback) to arbitrary linear functionals, 
albeit with different proofs.

\heading{References}

\def \nextref #1{\medskip \item {\bf [#1]}\quad }
 
\def\CA{Constr.\ Approx.}
\def\MC{Math.\ Comp.}
\def\MZ{Math.\ Z.}
 
\nextref {BR1}
C. de Boor and A. Ron\ (1990),
``On multivariate polynomial interpolation'',
{\it \CA \/} {\bf 6}, 287--302.
 
\nextref {BR2}
C. de Boor and A. Ron\ (1992),
``The least solution for the polynomial interpolation problem'',
{\it \MZ \/} {\bf 210}, 347--378.
 
\nextref {M}
C. A. Micchelli\ (1986),
``Interpolation of scattered data: distance matrices and conditionally positive definite functions'',
{\it \CA \/} {\bf 2}, 11--22.
 
\nextref {SX}
T. Sauer and Yuan Xu\ (1995),
``On multivariate Lagrange interpolation'',
{\it \MC \/} {\bf 64}, 1147--1170.
 
\nextref {S}
Robert Schaback\ (2002),
``Multivariate Interpolation by Polynomials and Radial Basis Functions'',
ms, April.
(available from Schaback's homepage)
\bye